\documentclass[reqno,11pt,a4paper]{amsart}

\usepackage{a4wide}

\usepackage{amsmath}
\usepackage{amssymb}
\usepackage{amsfonts}
\usepackage{amscd}
\usepackage{amsthm}
\usepackage{graphicx}
\usepackage[english]{babel}

\usepackage{url}
\usepackage{hyperref} 
\hypersetup{ 
  pdfauthor = {Romain Basson, Reynald Lercier, Christophe Ritzenthaler and Jeroen Sijsling}%
  pdftitle = {An Explicit Expression of the L\"uroth Invariant},%
  pdfsubject = {},%
  pdfkeywords = {Invariant ; ternary quartic ; genus 3 ; Ciani quartic ;
   algorithm ; Dixmier-Ohno invariants},%
  backref=true,%
  pagebackref=true,%
  hyperindex=true,%
  colorlinks=true,%
  breaklinks=true,%
  urlcolor=blue,%
  linkcolor=blue,%
  citecolor=blue,%
  bookmarks=true,%
  bookmarksopen=true}

\usepackage{breqn}
\newcounter{myequation}[equation]
\DeclareMathOperator{\GL}{GL}

\DeclareMathOperator{\SL}{SL}

\DeclareMathOperator{\Sym}{Sym}

\newcommand{\dd}{\mathcal{D}}

\newcommand{\ii}{\mathcal{I}}

\newcommand{\llm}{\mathcal{L}}

\newcommand{\qq}{\mathcal{Q}}

\newcommand{\FF}{\mathbb{F}}
\newcommand{\CC}{\mathbb{C}}

\newcommand{\PP}{\mathbb{P}}
\newcommand{\QQ}{\mathbb{Q}}

\newcommand{\ZZ}{\mathbb{Z}}

\title{An explicit expression of the L\"uroth invariant}

\author{Romain~Basson} 
\address{%
  IRMAR, %
  Universit\'e de Rennes 1, %
  Campus de Beaulieu, %
  35042 Rennes, %
  France. %
} 

\email{romain.basson@univ-rennes1.fr}

\author{Reynald Lercier}
\address{%
  \textsc{DGA MI}, %
  La Roche Marguerite, %
  35174 Bruz, %
  France. %
}

\address{%
  IRMAR, %
  Universit\'e de Rennes 1, %
  Campus de Beaulieu, %
  35042 Rennes, %
  France. %
}

\email{reynald.lercier@m4x.org}

\author{Christophe Ritzenthaler} 
\address{%
  IRMAR, %
  Universit\'e de Rennes 1, %
  Campus de Beaulieu, %
  35042 Rennes, %
  France. %
}
\email{ritzenth@iml.univ-mrs.fr}

\author{Jeroen~Sijsling} 
\address{Mathematics Institute, Zeeman Building, University of Warwick, %
Coventry CV4 7AL, United Kingdom.
}
\email{sijsling@gmail.com}

\thanks{The authors acknowledge support by the grants ANR-09-BLAN-0020-01 and PIEF-GA-2011-299887.}
\date{\today} \subjclass[2010]{14Q05 ; 13A50 ; 14H10 ; 14H37}
\keywords{invariant ; ternary quartic ; genus 3 ; Ciani quartic ;
   algorithm ; Dixmier-Ohno invariants}
\begin{document}

\begin{abstract}
  In this short note, we give an algorithm  to get an explicit expression of the L\"uroth invariant in terms of the Dixmier-Ohno invariants. We also get the explicit factorized expression on the locus of Ciani quartics in terms of the coefficients. Finally, we   answer two open questions on sub-loci of singular L\"uroth quartics. 
  
  \end{abstract}

\maketitle

\section{L\"uroth quartics}
\label{sec:introduction}

This note considers \emph{L\"uroth quartics}, which are plane quartics containing the ten vertices of a non-degenerate pentalateral. To make these notions precise, we give the following definitions. Let $V$ be a three-dimensional vector space over the complex field $\CC$.

\newtheorem{definition}{Definition}
\begin{definition}
A \emph{complete pentalateral} in the projective plane $\PP V$ is a curve $C \subset \PP V$ consisting of the union of five lines $\ell_1,\ldots,\ell_5$ that are three by three linearly independent (which is to say that the pairwise intersections of the lines $\ell_i$ yield exactly $10$ distinct points).

The \emph{vertices} of a complete pentalateral $C$ are the double points of $C$, that is, the $10$ points $\bigcup_{i \ne j} ((\ell_i=0) \cap (\ell_j=0))$. 
\end{definition}

\begin{definition}
Let $Q \subset \PP V$ be a non-singular quartic. Then $Q$ is called a \emph{non-singular L\"uroth quartic} if it contains the vertices of a complete pentalateral in $\PP V$.
\end{definition}

The set of plane quartics in $\PP V$ can be identified with the projective space $\PP \Sym^4 (V^*)$ over the fourth symmetric power of the dual vector space $V^*$ of $V$. This inherits an action of the group $\GL (V)$ and its derived subgroup $\SL (V)$, which act canonically on $V$. Choosing a basis, as we will do during our calculations, identifies $V$ with $\CC^3$ and the set of quartics $\PP \Sym^4 (V^*)$ with the projectivization of the vector space on the homogeneous degree $4$ monomials in the canonical basis $x$, $y$, $z$ of the dual space $(\CC^3)^*$. We will in turn identify this projective space with $\PP^{14}$ by choosing some ordering of these $15$ monomials. In this way, $\PP^{14}$ inherits an action of the groups $\GL_3 (\CC)$ and $\SL_3 (\CC)$.

The classical study of L\"uroth quartics culminated in 1919 with the work of Morley \cite{morley}. This showed that the Zariski closure of the locus of non-singular L\"uroth quartics in the projective space $\PP \Sym^4 (V^*)$ is an irreducible hypersurface described as the vanishing locus of a single homogeneous polynomial function $L$ on the projective space of quartics $\PP \Sym^4 (V^*)$, well-defined up to scalars. We shall call this $L$ the \emph{L\"uroth invariant}. Morley showed that $L$ is of degree $54$.

\begin{definition}
Let $Q \subset \PP V$ be a quartic. Then $Q$ is called a \emph{L\"uroth quartic} if $L(Q) = 0$.
\end{definition}
In recent years, after the seminal work of \cite{barth}, several authors have revived this subject in~\cite{boehning,ottaviani2,ottaviani,ottaviani3} (see also \cite{shakirov} on the undulation invariant).

However, an explicit expression of $L$ was still missing. In the following section we explain how to compute such an expression. Our main new technique lies in an effective use of \cite{ohno}, an unfortunately unpublished article in which Ohno gives a complete set of generators for the invariants of ternary quartics under the action of $\SL_3(\CC)$, completing the set of primary invariants found in \cite{dixmier}.
These invariants were also used in~\cite{giko}, and new effective methods to verify their correctness can be found in~\cite{elsenhaus}. Our calculations use the implementation of these invariants in \texttt{Magma}, which is due to Kohel.

\section{The algorithm}
\label{sec:mainalgo}

The key point is the following observation:
\newtheorem{proposition}{Proposition}
\begin{proposition}\label{prop:LInv}
The homogeneous polynomial function $L$ on $\PP \Sym^4 (V^*)$ is $\GL(V)$-invariant up to scalars. In particular, $L$ is $\SL (V)$-invariant.
\end{proposition}
\begin{proof}
Since we are working over an algebraically closed field, this is obvious from the obvious fact that any $\GL(V)$-transform of a L\"uroth quartic is again a L\"uroth quartic. 
\end{proof}

Let 
\begin{align*}
R = S [ \PP \Sym^4 (V^*) ]^{\SL (V)}
\end{align*}
be the ring of $\SL(V)$-invariant homogeneous polynomial functions on $\PP \Sym^4 (V^*)$, which coincides with those functions that are $\GL(V)$-invariant up to scalars. The structure of the ring $R$ is known. Indeed, let $I = (I_3, I_6, I_9,  I_{12}, I_{15}, I_{18},  I_{27})$ be the primary invariants for ternary quartics under the action of $\SL_3(\CC)$ found by Dixmier in~\cite{dixmier}, and let $J = (J_9,J_{12},  J_{15},  J_{18}, I_{21},J_{21})$ be the secondary invariants found by Ohno~\cite{ohno}. In both cases, the index specifies the degree of the invariant as a homogeneous function. Then we have:

\newtheorem{theorem}{Theorem}
\begin{theorem}[Dixmier-Ohno]\label{thm:RInv} We have $R = \CC [I,J]$.
\end{theorem}

We now choose a basis for $V$ and a corresponding coordinatization of $\PP \Sym^4 (V^*) \cong \PP^{14}$ as the projective space over the degree $4$ monomials in $x$, $y$, $z$. The function $L$ then becomes a homogeneous expression in the coefficients of these monomials. It is unlikely that this function can be written down in any reasonable way (see also the final remark in Section~\ref{sec:open}). However, by Proposition~\ref{prop:LInv} and Theorem~\ref{thm:RInv}, we can express $L$ as a polynomial in the invariant functions in $I$ and $J$. This expression is not unique, since as we shall see, there are relations between these invariant monomials in degree $54$.

To obtain an expression for $L$, we apply the method of evaluation-interpolation. This is based on the following observation:
\begin{proposition}
Let $S = \CC [x_3 , \ldots , x_{27}, y_9 , \ldots , y_{21} , y'_{21}]$ be a graded polynomial algebra in $13$ variables, weighted by indices, and consider the surjection $q$ given by
\begin{align*}
q : S & \longrightarrow R  \\
x_k & \longmapsto I_k \\
y_{\ell} & \longmapsto J_{\ell}  .
\end{align*}
Let $R_{54} \subset R$ be the set of homogeneous functions of degree $54$, and define $S_{54} \subset S$ analogously. Let $K$ be the kernel of the map $S_{54} \rightarrow R_{54}$. Then the $\mathrm{dim} (K) = 215$.

Let $X$ be a finite set of L\"uroth quartics. Consider the linear map $q' : S_{54} \rightarrow \CC^X$ given by evaluating at the polynomials in $X$. Let $K'$ be the kernel of $q'$, and suppose that $\mathrm{dim} (K') = 216$. Let $\widetilde{L}$ be an element of $K' \backslash K$. Then the image $L = q(\widetilde{L})$ equals the L\"uroth invariant.
\end{proposition}
\begin{proof}
One calculates that the $\mathrm{dim} (S_{54} ) = 1380$. Since calculating the Hilbert polynomial of $R$ as in~\cite[p.1045]{shioda67} yields $\mathrm{dim} (S_{54} ) = 1165$, we indeed find that $\mathrm{dim} (K) = 215$. The rest is straightforward considering the uniqueness of $L$ up to scalars.
\end{proof}

The details of the calculation are therefore as follows.
\begin{enumerate}
\item Construct the $1380$ monomials $$\ii=\{I_3^{18},I_3^{16} I_6, I_3^{15} I_9, I_3^{15} J_9, \ldots,J_{18}^3, I_{27}^2\}$$ of degree $54$ that generate the $\CC$-vector space of invariants of degree $54$.
\item Generate a sufficiently large finite set $\qq$ of cardinality $q$ of random plane quartics with rational coefficients.
\item Generate a sufficiently large finite set $\llm$ of cardinality $l$ of random L\"uroth quartics of the form 
\begin{dgroup*}[compact,style={\small},spread={0pt}]
  \begin{dmath*}
\ell_1 \ell_2 \ell_3 \ell_4 + c_1 \cdot \ell_2 \ell_3 \ell _4 \ell_5 + c_2 \cdot \ell_1 \ell_3 \ell_4 \ell_5 + c_3 \cdot \ell_1 \ell_2 \ell_4 \ell_5 + c_4 \cdot  \ell_1 \ell_2 \ell_3 \ell_5
\end{dmath*}
\end{dgroup*}
where $\ell_1=x,\ell_2=y,\ell_3=z,\ell_4=x+y+z$, $\ell_5$ is a line with random rational coefficients and $c_i$ are rational coefficients.
\item Compute the matrix $M_1= (I(q))_{I \in \ii,q \in \qq}$, evaluating the monomials in $\ii$ at the quartics in $\qq$.
\item Compute the matrix $M_2= (I(q))_{I \in \ii,q \in \llm}$, evaluating the monomials in $\ii$ at the L\"uroth quartics in $\llm$.
\item Compute the $215$-dimensional kernel $N_1$ of $M_1$. This gives a basis of the homogeneous relations of degree $54$ that are satisfied by the invariants of all ternary quartics.
\item Compute the $216$-dimensional kernel $N_2$ of $M_2$. This gives a basis of the homogeneous relations of degree $54$ that are satisfied by all L\"uroth quartics.
\item A non-zero element in the complement of $N_2$ in $N_1$ is an expression for $L$ in terms of the Dixmier-Ohno invariants.
\end{enumerate}
All these computations were done with Magma software. Over finite fields
$\FF_p$ with prime cardinality $p= 2 017$, $10 007$, $100 003$ or even $1 000 003$, computations can
be done in less than a minute. However, getting the result over the rationals
is more challenging. The main concern is to deal with
matrices $M_1$ and $M_2$ whose coefficients are as small as possible. So, at Step
(2) of the algorithm, we generate plane quartics with random integer
coefficients only equal to $-1$, $0$, or $1$. Similarly, we restrict Step (3)
to L\"uroth forms defined by integer coefficients $c_i$ bounded in absolute value
by 4.

We can estimate the size of the computations involved in this run of the algorithm by using the
Hadamard bounds for our matrices $M_1$ and $M_2$; the quartics under consideration
yield bounds slightly smaller than $2^{200\,000}$ for $M_1$ and $2^{350\,
  000}$ for $M_2$. As a sanity check before running the code over the rationals,
we verify that this subset of quartics yields a valid result modulo small primes.

Most of the time is spent at Step (6) and Step (7) of the
algorithm, precisely 5 and 9 hours on our laptop (based on a \textsc{Intel Core
  i7 M620 2.67GHz} processor).

A program to get the result is available on the web page of the authors\footnote{
\small \url{http://perso.univ-rennes1.fr/christophe.ritzenthaler/programme/luroth/luroth.m}
}.
It uses the implementation of the Dixmier-Ohno invariants in Magma by Kohel\footnote{
\small \url{http://echidna.maths.usyd.edu.au/kohel/alg/index.html}}.
The $1.4$Mb result is also available online\footnote{
 \small \url{http://perso.univ-rennes1.fr/christophe.ritzenthaler/programme/luroth/LurothInvF.m}}.
 It is given by $1164$ monomials with rational coefficients, the largest of which is a quotient of a $680$-digit integer by a coprime $671$-digit integer. Modulo $1000003$, the expression starts as
\begin{dgroup*}[compact,style={\small},spread={0pt}]
  \begin{dmath*}
I_3^{18} + 469313 I_3^2 I_6^8 + 710780 I_6^9 + 969230 I_3^3 I_6^6 I_9 + 374233 I_3 I_6^7 I_9
    + 276144 I_3^2 I_6^5 I_9^2 + 602674 I_6^6 I_9^2 + 527614 I_3^3 I_6^3 I_9^3 + 
    538637 I_3 I_6^4 I_9^3 + 392526 I_3^4 I_6 I_9^4 + 645841 I_3^2 I_6^2 I_9^4 + 
    914224 I_6^3 I_9^4 + 207808 I_3^3 I_9^5 + 31577 I_3 I_6 I_9^5 + 635768 I_9^6 + 
    668878 I_3^{15} J_9 + 507293 I_3^3 I_6^6 J_9 + 318476 I_3 I_6^7 J_9 + 
    59775 I_3^2 I_6^5 I_9 J_9 + 581086 I_6^6 I_9 J_9 + 830307 I_3^3 I_6^3 I_9^2 J_9 + 
    804817 I_3 I_6^4 I_9^2 J_9 + 6418 I_3^6 I_9^3 J_9 + 578316 I_3^4 I_6 I_9^3 J_9 + 
    741618 I_3^2 I_6^2 I_9^3 J_9 + 452974 I_6^3 I_9^3 J_9 + 36214 I_3^3 I_9^4 J_9 + 
    522408 I_3 I_6 I_9^4 J_9 + 253043 I_9^5 J_9 + 469299 I_3^2 I_6^5 J_9^2 +\ldots 
\end{dmath*}
\end{dgroup*}

\section{Ciani quartics}
We call a \emph{Ciani quartic} a plane quartic of the form 
$$a x^4 +  b x^2  y^2 + c x^2 z^2 + d y^4 + e y^2 z^2 + f z^4.$$ 
A generic Ciani quartic has automorphism group isomorphic with $\ZZ / 2 \ZZ \times \ZZ / 2 \ZZ$, and conversely every quartic with this property is $\CC$-isomorphic to a Ciani quartic. The dimension of the substratum of Ciani quartics in the full dimension $6$ moduli space of plane quartics equals $3$.

In \cite[Sec.5]{hauenstein}, using  different techniques,  Hauenstein and Sottile obtained the factorization on Ciani quartics  of the L\"uroth  invariant as
$$G^4   H^2  J$$
with $G,H,J \in \CC[a,b,c,d,e,f] $ homogeneous of respective degree $6,9$ and $12$. Using our expression, it is easy to confirm their decomposition. We give a slightly different version of the result, which is due to the fact the coefficients $b,c,e$ are replaced by $2b,2c,2e$ in \cite[Sec.5]{hauenstein}: 
  $$G=a \cdot  d \cdot f \cdot (a d f - (1/4) a e^2 - (1/4) b^2 f - (1/4) b c e - (1/4) c^2 d),$$
\begin{dgroup*}[compact,style={\small},spread={0pt}]
  \begin{dmath*}
 H= (a d f - (1/4) a e^2 - (1/4) b^2 f + (1/4) b c e + (3/4) c^2 d) \cdot
  (a d f - (1/4) a e^2 + (3/4) b^2 f + (1/4) b c e - (1/4) c^2 d) \cdot
  (a d f + (3/4) a e^2 - (1/4) b^2 f + (1/4) b c e - (1/4) c^2 d),
\end{dmath*}
\end{dgroup*}
\begin{dgroup*}[compact,style={\small},spread={0pt}]
  \begin{dmath*}
J=  a^4 d^4 f^4 - (1/49) a^4 d^3 e^2 f^3 + (51/19208) a^4 d^2 e^4 f^2 -
      (1/38416) a^4 d e^6 f + (1/614656) a^4 e^8 - (1/49) a^3 b^2 d^3 f^4 -
      (205/9604) a^3 b^2 d^2 e^2 f^3 - (3/38416) a^3 b^2 d e^4 f^2 +
      (1/153664) a^3 b^2 e^6 f + (15/343) a^3 b c d^3 e f^3 +
      (29/9604) a^3 b c d^2 e^3 f^2 - (5/38416) a^3 b c d e^5 f -
      (1/153664) a^3 b c e^7 - (1/49) a^3 c^2 d^4 f^3 -
      (205/9604) a^3 c^2 d^3 e^2 f^2 - (3/38416) a^3 c^2 d^2 e^4 f +
      (1/153664) a^3 c^2 d e^6 + (51/19208) a^2 b^4 d^2 f^4 -
      (3/38416) a^2 b^4 d e^2 f^3 + (3/307328) a^2 b^4 e^4 f^2 +
      (29/9604) a^2 b^3 c d^2 e f^3 - (5/19208) a^2 b^3 c d e^3 f^2 -
      (3/153664) a^2 b^3 c e^5 f - (205/9604) a^2 b^2 c^2 d^3 f^3 +
      (2/2401) a^2 b^2 c^2 d^2 e^2 f^2 + (55/153664) a^2 b^2 c^2 d e^4 f +
      (3/307328) a^2 b^2 c^2 e^6 + (29/9604) a^2 b c^3 d^3 e f^2 -
      (5/19208) a^2 b c^3 d^2 e^3 f - (3/153664) a^2 b c^3 d e^5 +
      (51/19208) a^2 c^4 d^4 f^2 - (3/38416) a^2 c^4 d^3 e^2 f +
      (3/307328) a^2 c^4 d^2 e^4 - (1/38416) a b^6 d f^4 +  (1/153664) a b^6 e^2 f^3
      - (5/38416) a b^5 c d e f^3 - (3/153664) a b^5 c e^3 f^2 -
      (3/38416) a b^4 c^2 d^2 f^3 + (55/153664) a b^4 c^2 d e^2 f^2 +
      (3/153664) a b^4 c^2 e^4 f - (5/19208) a b^3 c^3 d^2 e f^2 -
      (17/76832) a b^3 c^3 d e^3 f - (1/153664) a b^3 c^3 e^5 -
      (3/38416) a b^2 c^4 d^3 f^2 + (55/153664) a b^2 c^4 d^2 e^2 f +
      (3/153664) a b^2 c^4 d e^4 - (5/38416) a b c^5 d^3 e f -
      (3/153664) a b c^5 d^2 e^3 - (1/38416) a c^6 d^4 f +  (1/153664) a c^6 d^3 e^2
      + (1/614656) b^8 f^4 - (1/153664) b^7 c e f^3 + (1/153664) b^6 c^2 d f^3 +
      (3/307328) b^6 c^2 e^2 f^2 - (3/153664) b^5 c^3 d e f^2 -
      (1/153664) b^5 c^3 e^3 f + (3/307328) b^4 c^4 d^2 f^2 +
      (3/153664) b^4 c^4 d e^2 f + (1/614656) b^4 c^4 e^4 -
      (3/153664) b^3 c^5 d^2 e f - (1/153664) b^3 c^5 d e^3 +
      (1/153664) b^2 c^6 d^3 f + (3/307328) b^2 c^6 d^2 e^2 -  (1/153664) b c^7 d^3 e
      + (1/614656) c^8 d^4.
\end{dmath*}
\end{dgroup*}
The product $G^4 H^2 J$ has $1695$ monomials. Note that the total amount of  weighted monomials in $a,b,c,d,e$ and $f$ in a generic degree $54$ invariant is $3439$.

\section{Singular L{\"u}roth quartics}
\label{sec:singular}

Let $\llm \subset \PP \Sym^4 (V^*)$ be the locus of L\"uroth quartics, and let $\dd \subset \PP \Sym^4 (V^*)$ be the discriminantal hypersurface defined by the equation $I_{27} = 0$. We will now obtain new results on the geometry of the locus $\llm \cap \dd$ of \emph{singular L\"uroth quartics}. Work by Le Potier and Tikhomirov \cite{LePT} shows that
\begin{align*}
\llm \cap \dd = \llm_1 \cup \llm_2  ,
\end{align*}
where $\llm_1$ and $\llm_2$ are irreducible subschemes of $\PP \Sym^4 (V^*)$ of codimension $2$ whose respective degrees as subschemes of $\dd$ equal $24$ and $30$ respectively. Moreover, while $\llm_1$ is reduced, the reduced subscheme $(\llm_2)_{\mathrm{red}}$ of $\llm_2$ is of degree $15$.

In \cite{ottaviani}, Ottaviani and Sernesi showed that no new degree $15$ invariant vanishes on $(\llm_2)_{\mathrm{red}}$, which implies that this scheme is not a principal hypersurface in $\llm$. We will prove a stronger result, namely that none of $\llm_1$, $\llm_2$, $(\llm_2)_{\mathrm{red}}$ is a complete intersection. We apply the same methods as in Section~\ref{sec:mainalgo}. The main problem now is to generate quartics in $\llm_1$ and $\llm_2$. 

For $\llm_2$, we can proceed by using Remark 3.3 in \cite{ottaviani}: we now choose the lines in Step (3) above such that three of them have a common point of intersection.

For $\llm_1$, we have to perform the constructions and results in \cite[p. 1759]{ottaviani}. The procedure is as follows.
\begin{enumerate}
\item Construct a cubic surface $S$ with two skew lines $l$, $m$.
\item Calculate the double cover $f : l \rightarrow m$ sending $p \in l$ to the intersection $T_p S \cap m$ of the tangent plane $T_p S$ to $p$ at $S$ with the line $m$, and construct $g : m \rightarrow l$ analogously.
\item Let $B_f \subset m$ (resp. $B_g \subset l$) be the branch divisor of $f$ (resp. $g$). Construct morphisms $f' : l \rightarrow \PP^1$ (resp. $g' : m \rightarrow \PP^1$) ramifying over $B_g$ (resp. $B_f$).
\item Construct $Q \in m \subset S$ such that $f^{-1} (Q)$ is also a fiber of $f'$.
\item Construct the ramification locus of the degree $2$ projection $S \rightarrow \PP (T_Q S)$ from the point $Q \in S$. Then by Proposition 3.1(i) of \cite{ottaviani}, we obtain a quartic in $\llm_1$.
\end{enumerate}

The following string of propositions and remarks show how these steps can be implemented without extending the base field (we will take $k = \QQ$ throughout).

\begin{proposition}
Let $k$ be a field, and suppose that we are given six rational points $p_1 , \ldots , p_6 \in \PP V (k)$ in the projective plane over $k$. Suppose additionally that this set of points is sufficiently general in the sense that the complete linear system $C$ of cubics passing through them has dimension $4$. Construct the Clebsch rational map $c : \PP V \rightarrow \PP C$, and let $S \subset \PP C$ be the Zariski closure of $c(\PP V)$. Then $S$ is the vanishing locus of a quaternary cubic form $F \in \Sym^3 (C^*)$ over $k$.

The rational map $c$ restricts to a birational map between $\PP V$ and $S$. Let $l_0 \subset \PP V$ be the rational line containing $p_1$ and $p_2$, and let $m_0 \subset \PP V$ be the rational line containing $p_1$ and $p_3$. Let $x$ and $x'$ be two points in $l$ not equal to $p_1$ or $p_2$, and let $y$ and $y'$ be two points in $l$ not equal to $p_1$ or $p_3$. Then the images $c(x)$ and $c(x')$ are well-defined elements of $S$, and the line $l$ through them is defined over $k$ and included in $S$. Analogously, one obtains a line $m$ through $c(y)$ and $c(y')$. The lines $l$ and $m$ are skew.
\end{proposition}
\begin{proof}
This is a standard result from the theory of cubic surfaces, see~\cite[Section V.4]{hartshorne}. 
\end{proof}

This deals with part 1. To perform the calculations in point 2 explicitly, choose coordinates on $l$ by taking two points $l_1,l_2$ on $l$ and sending $(x:y) \in \PP^1$ to $x l_1 + y l_2$, and similarly on $m$ by choosing $m_1,m_2 \in m$. To determine the morphism $f : l \rightarrow m$ explicitly in these in coordinates, we choose two equations $M_1 = M_2 = 0$ defining $m$. Given $p \in l$ with coordinates $(x:y) \in \PP^1$, the point $f(p) = T_p S \cap m$ corresponds to the vector space that is the kernel of the matrix whose rows are given by $M_1$, $M_2$ and the partial derivatives of $F$. A generating vector for this space will be a combination of $m_1$ and $m_2$ with homogeneous quadratic coefficients $f_1(x,y),f_2(x,y)$ in $(x,y)$. The morphism $f$ now corresponds to the map $\PP^1 \rightarrow \PP^1$ given by $(f_1,f_2)$. Similarly, one determines $g$. For point $3$, we use the following result.

\begin{proposition}\label{prop:invo}
Let $l$ be a projective line over $k$, with homogeneous coordinates $x$ and $y$, and let $D$ be a $k$-rational divisor of degree $2$ on $l$. Then the following formulae determine a degree $2$ morphism $f' : l \rightarrow \PP^1$ over $k$ whose ramification locus equals $D$.
\begin{itemize}
\item If $D$ consists of two points $(x_1 : y_1)$ and $(x_2 : y_2)$ that are rational over $k$, then one can take
\begin{align*}
f'(x:y) = ((y_1 x - x_1 y)^2:(y_2 x - x_2 y)^2) .
\end{align*}
\item If $D$ is defined by an equation $r x^2 + t y^2 = 0$, then one can take
\begin{align*}
f'(x:y) = (r x^2 - 2 t x y - t y^2 : r x^2 + 2 t x y - t y^2) ;
\end{align*}
\item If $D$ is defined by an equation $r x^2 + s x y +  t y^2 = 0$ with $s \neq 0$, then one can take 
\begin{align*}
f'(x:y) = & (r^2 s x^2 + 2 r (s^2 - 2 r t) x y + (s^3 - 3 r s t) y^2 : \\ & \; r ( r s x^2 + 4 r t x y + s t y^2)) .
\end{align*}
\end{itemize}
\end{proposition}
\begin{proof}
Once the answer is given, the verification is trivial. But let us illustrate how to find these expressions by treating the thir case, where the points of $D$ are defined over a proper quadratic extension of $k$. We use the affine coordinate $t = x/y$ on $l$. Suppose that using this coordinate, the divisor $D = [d] + [\overline{d}]$ consists of two conjugate points, not summing to zero because we are in case (iii). Then $f' = ((t-d)/(t-\overline{d}))^2$ of above now has the property that $\overline{f'} = 1 / f'$. But then one verifies that $(d f' + \overline{d})/(\overline{d} f' + d)$ is a fractional linear transformation of $f'$ that is stable under conjugation and hence defines a morphism over the ground field. Homogenizing, one obtains the more elegant expression given in the statement of the proposition.
\end{proof}

We now treat point 4.

\begin{proposition}
Let $f,f' : \PP^1 \rightarrow \PP^1$ be two degree $2$ morphisms, neither of which can be obtained from the other by postcomposing with an automorphism of $\PP^1$. Then the $q$ in $\PP^1$ such that the fiber of $f$ over $q$ is also a fiber of $f'$ over a point $q'$ can be obtained as follows.

Write $$f (x:y) = (a_1 x^2 + b_1 x y + c_1 y^2 : a_2 x^2 + b_2 x y + c_2 y^2)$$ and $$f' (x:y) = (a'_1 x^2 + b'_1 x y + c'_1 y^2 : a'_2 x^2 + b'_2 x y + c'_2 y^2).$$ Then $q = (\lambda_1 : \lambda_2)$, where $(\lambda_1 , \lambda_2 , \lambda'_1 , \lambda'_2)$ generates the kernel of the matrix
\begin{align*}
\left(\begin{array}{cccc} a_1 & -a_2 & -a'_1 & a'_2 \\ b_1 & -b_2 & -b'_1 & b'_2 \\ c_1 & -c_2 & -c'_1 & c'_2 \end{array}\right).
\end{align*}
\end{proposition}
\begin{proof}
If we let $q = (\lambda_1 : \lambda_2)$ and $q' = (\lambda'_1 : \lambda'_2)$, then finding $q$ and $q'$ comes down to solving the equation
\begin{align*}
 \lambda_2 (a_1 x^2 + b_1 x y + c_1 y^2)   - \lambda_1 (a_2 x^2 + b_2 x y + c_2 y^2) \\
= \lambda'_2 (a'_1 x^2 + b'_1 x y + c'_1 y^2) - \lambda'_1 (a'_2 x^2 + b'_2 x y + c'_2 y^2),
\end{align*}
which evidently corresponds to the determination of the kernel of the matrix in question.
\end{proof}

To calculate point 5, we choose an isomorphism $\PP C \cong \PP^3$ mapping the point $Q$ to $(1:0:0:0)$ and apply the following elementary result.

\begin{proposition} 
Let $S \subset \PP^3$ be a cubic surface containing $Q = (1:0:0:0)$ that is defined by a quaternary cubic form $F \in k[w,x,y,z]$. Let $q$ be the projection $S \rightarrow \PP (T_Q S)$ from the point $Q \in S$. Then the ramification locus of $q$ is isomorphic to the quartic curve in $\PP^2$ determined by the vanishing of the discriminant of the quadratic polynomial $F(1,xt,yt,zt)/t$.
\end{proposition}
\begin{proof}
Since $S$ is a subscheme of $\PP^3$, we get an induced coordinatization of $T_Q S$ by sending $(x:y:z)$ to the tangent direction given by the line through the points $(1:0:0:0)$ and $(1:x:y:z)$. Then $(x:y:z)$ is a ramification point of the projection $S \rightarrow \PP (T_Q S)$ if and only if the equation $F(1,xt,yt,zt) = 0$ has a double root outside $0$, or in other words if the discriminant of the quadratic polynomial $F(1,xt,yt,zt)/t$ vanishes. This discriminant is a homogeneous quartic form in the variables $x,y,z$, which defines the plane quartic in $\llm_1$ that we were looking for.
\end{proof}

A program to generate quartics in $\llm_1$ by using the steps above is available online\footnote{\small \url{http://perso.univ-rennes1.fr/christophe.ritzenthaler/programme/luroth/GenerateL1.m}}.

\newtheorem{remark}{Remark}
\begin{remark}We also tried to generate quartics in $\llm_1$ by using Remark 3.4 of \cite{ottaviani}. We take cubics $S$ of the form
$$t^2 x + t (a x^2 + 2 b  xy + 2 c xz+ d y^2+2  e yz+f z^2) + g(x,y,z)$$
with $g$ a random degree $3$ homogeneous polynomial, such that $S$ is non singular and $e^2=d f$. The last condition ensures that  $p=(0:0:0:1)$ belongs to the Hessian $H$ of $S$ and, after checking that $p$ is non singular, we take quartics which are tangent plane sections of $H$ at $p$. Unfortunately, it seems that these quartics are special in $L_1$, since there are degree $24$ relations between their invariants (there is a $27$ dimensional space of relations in degree $24$ between randomly generated quartics of this form).
\end{remark}

Having generated a sufficiently large database\footnote{$1024$ curves over $\QQ$ are available at 
\url{http://perso.univ-rennes1.fr/christophe.ritzenthaler/programme/luroth/L1Database.m}} of curves in $\llm _1$ by choosing random $6$-tuples $\{p_1 , \ldots , p_6\}$, we can again proceed as in Section \ref{sec:mainalgo}. Up to degree $30$, all invariants vanishing on the quartics in this databases for $\llm _1$ and $\llm _2$ are multiples of $I_{27}$. Since the codimension $2$ components $\llm _1$, $\llm _2$, $(\llm_2)_{\mathrm{red}}$ of $\llm \cap \dd$ have degree at most $24 \cdot 27$,$24 \cdot 30$,$24 \cdot 15$, which are all smaller than $(30)^2$, we have the following result.

\begin{theorem}\label{thm:IntThm}
The subschemes $\llm _1$, $\llm _2$, $(\llm_2)_{\mathrm{red}}$ of the projective space of quartic curves $\PP \Sym^4 (V^*)$ (and hence their images in the coarse moduli space of plane quartic curves) are not complete intersections. In particular, they are not principal hypersurfaces in the discriminant locus $\dd$.
\end{theorem}

As there is no degree $24$ invariant vanishing on $\llm_1$, Morley's putative construction of such an invariant $I_{24}$ in \cite[p.282]{morley} is incorrect. On the authors' webpage,
a Magma program\footnote{\small \url{http://perso.univ-rennes1.fr/christophe.ritzenthaler/programme/luroth/SingularLurothInv.m}} is available to check all steps on the way to Theorem~\ref{thm:IntThm}.

\section{Open questions}\label{sec:open}
The  expression $L$ of the L\"uroth invariant that we found depends on several arbitrary choices that may explain its cumbersomeness. First, there is the choice of the basis of invariants. Though some of the Dixmier invariants have geometrical interpretations that are `natural', the same is far from evident for the new Ohno invariants. Secondly, our choice can be modified by any element of the kernel $N_1$. Beyond a cancellation of the coefficients of $215$ of these monomials that we have already accomplished by simple linear algebra, further minimization of the number of monomials in the expression for $L$ could in theory be achieved by techniques based on coding theory. Still, the parameters seem too large to make this feasible in practice.\\
The negative answers concerning the existence of degree 24 and 30 invariants in Section \ref{sec:singular} exclude the decomposition from \cite[p.1764]{ottaviani}. The geometry of the situation does not seem to give a clue for the existence of another such decomposition. \\
An expression in terms of the $15$ coefficients of the generic quartic would of course be useful. However, it is not even practically achievable to formally express the fundamental Dixmier-Ohno invariants in this way, since these expressions contain far too many monomials, as is for instance the case for the discriminant $I_{27}$. A count of weighted monomials in $15$ variables for degree $54$ invariants leads to a total of $62 \, 422 \, 531 \, 333$. Of course only a fraction of these monomials may occur in the final expression of $L$, but we could not figure out their number, let alone the Newton polytope of $L$.

\par\medskip
\noindent {\bf Acknowledgments.} This note arose from Giorgio Ottaviani's great lectures at a workshop on invariant theory and projective geometry held in Trento in September 2012, in which the main problems treated in it were mentioned. The authors wish cordially to thank Ottaviani for his further elaborations on the subject, as well as his careful reading of a first draft of this article.

The authors acknowledge support by the grants ANR-09-BLAN-0020-01. The fourth author was additionally supported by the Marie Curie Fellowship IEF-GA-2011-299887.

\end{document}